\documentclass[11pt,a4paper]{article}
\usepackage[T1]{fontenc}
\usepackage{lmodern}
\usepackage[margin=29mm]{geometry}
\usepackage{amsmath,amssymb,amsthm,mathtools}
\usepackage{microtype}
\usepackage{xurl}
\usepackage[hidelinks]{hyperref}
\hypersetup{pdftitle={Every Cubic Graph Admits a 2-Quasi-Bisection},pdfauthor={Aleksandr Mamleev and Daniil Shapko},pdfsubject={Graph theory and quasi-bisections},pdfkeywords={cubic graphs, cactus, graph coloring, Lean}}
\usepackage{booktabs,tabularx}
\newtheorem{theorem}{Theorem}
\newtheorem{lemma}[theorem]{Lemma}
\newtheorem{corollary}[theorem]{Corollary}
\theoremstyle{definition}

\theoremstyle{remark}
\newtheorem{remark}[theorem]{Remark}
\newcommand{\sgnset}{\{-1,+1\}}
\newcommand{\ind}[1]{\mathbf{1}_{\{#1\}}}
\newcommand{\abs}[1]{\lvert #1\rvert}
\newcommand{\N}{\mathbb{N}}
\newcommand{\Z}{\mathbb{Z}}
\newcommand{\lean}[1]{\nolinkurl{#1}}
\title{Every Cubic Graph Admits a 2-Quasi-Bisection}
\author{Aleksandr Mamleev\textsuperscript{1} and Daniil Shapko\textsuperscript{2}\\[5pt]
\small \textsuperscript{1}National Research Nuclear University MEPhI\\
\small (Moscow Engineering Physics Institute), Moscow, Russia\\[2pt]
\small \textsuperscript{2}MIREA -- Russian Technological University, Moscow, Russia}
\date{}
\begin{document}
\maketitle
\begin{abstract}
We prove the 2-quasi-bisection conjecture of Esperet, Mazzuoccolo and Tarsi for finite simple cubic graphs. More generally, every finite simple subcubic graph of even order has a two-colouring whose colour classes differ in size by at most two and whose monochromatic components have at most two vertices. The principal ingredient is a sharp bound on the gaps between attainable colour imbalances of a subcubic cactus with prescribed boundary colours. We establish a weighted version of this bound, in which the weights are coupled to the boundary conditions. This coupling permits contraction across bridges without increasing the gap bound. A maximal packing of vertex-disjoint even cycles reduces the general problem to the cactus case. By inverting these cycles one vertex at a time and tracking the imbalance on the uncovered vertices, we either obtain a suitable colouring directly or complete it by repair chains starting at the endpoints of at most two marked edges. The proof is independent of finite computational searches. The accompanying Lean formalization proves the even-order subcubic theorem and derives the cubic case as a corollary.
\end{abstract}

\section{Introduction}
All graphs are finite, simple and undirected unless stated otherwise. A graph is subcubic if its maximum degree is at most three, and cubic if every vertex has degree three. For a colouring $h:V(G)\to\sgnset$, define
\begin{equation}\label{eq:ordinary}
 D_G(h)=\sum_{v\in V(G)}h(v),\qquad
 d_h(v)=\abs{\{u\in N_G(v):h(u)=h(v)\}}.
\end{equation}
We call $h$ admissible if $d_h(v)\leq 1$ for every vertex. In a simple graph this is equivalent to requiring that every monochromatic component has at most two vertices: a connected graph on at least three vertices has a vertex of degree at least two. A 2-quasi-bisection is an admissible colouring with $\abs{D_G(h)}\leq 2$.

Esperet, Mazzuoccolo and Tarsi conjectured that every cubic graph admits a 2-quasi-bisection \cite[Conjecture 10]{EMT}. Liu, Hao, Luo and Zhang proved the conjecture for bridgeless cubic graphs of oddness two \cite{LHLZ}. DeVos and Nurse \cite{DN} and Tang, Hao, Luo and Zhang \cite{THLZ} obtained partial results under additional decomposition hypotheses. Our proof uses local defect repairs related to the switching arguments in \cite[Section 4]{THLZ}, but relies on a weighted bound on gaps between attainable colour imbalances of cacti to remove the decomposition restrictions.

\begin{theorem}[Main theorem]\label{thm:main}
Every finite simple cubic graph $G$ admits a colouring $h:V(G)\to\sgnset$ such that
\[
 \abs{D_G(h)}\leq 2\quad\text{and}\quad d_h(v)\leq1\quad(v\in V(G)).
\]
No connectedness or bridgelessness assumption is required.
\end{theorem}
We will in fact prove the following extension, using exactly the same boundary model.
\begin{theorem}\label{thm:subcubic}
Every finite simple subcubic graph of even order admits a 2-quasi-bisection.
\end{theorem}
The argument does not assert the existence of an exact bisection, for which the two classes must have equal size.

The proof proceeds as follows. Choose a maximal family of vertex-disjoint even cycles and colour each of them alternately. The uncovered vertices induce a disjoint union of cacti. We invert the chosen cycles one vertex at a time and track the imbalance on the uncovered vertices, restoring admissibility there after each flip. This imbalance is even and changes by at most two at each such step. If it reaches zero during an inversion, any defects are confined to the endpoints of the two boundary edges of the inverted segment, and repair chains yield a 2-quasi-bisection. Otherwise, after all cycles have been inverted, a gap bound of four for the attainable imbalances of the residue allows us to pass to the inverse of the initial colouring of the residue without skipping all of $-2,0,2$. The main technical task is to prove this gap bound uniformly for cacti with prescribed boundary colours.

\section{Boundary models and defect repair}\label{sec:boundary}
A boundary model consists of a finite simple graph $H$ of maximum degree at most three, together with three nonnegative integer parameters at each vertex: $f_v(+1),f_v(-1)$, and $o_v$. We also specify a number $m_v$ of marked opposite positions. These data satisfy
\begin{equation}\label{eq:budget}
 \deg_H(v)+f_v(+1)+f_v(-1)+o_v=3,\qquad 0\leq m_v\leq o_v.
\end{equation}
The $f_v(t)$ positions have a fixed colour $t$. The $o_v$ opposite positions are always assigned the colour $-h(v)$, so never contribute a monochromatic incidence. Exactly $m_v$ of these positions are marked. Equivalently, each opposite position carries a label in $\{0,1\}$ and $m_v$ is the sum of these labels. These positions are local data, not new vertices or edges of $H$.

Write
\begin{equation}\label{eq:model}
 \begin{aligned}
 d_h(v)&=\sum_{u\in N_H(v)}\ind{h(u)=h(v)}+f_v(h(v)),\\
 w_v&=1-m_v,\qquad S(h)=\sum_{v\in V(H)}w_vh(v).
 \end{aligned}
\end{equation}
Admissibility means $d_h(v)\leq1$ at every vertex, now including the fixed boundary incidences. A vertex with $d_h(v)\geq2$ is a defect. The possible weights are $-2,-1,0,1$. Negative weights are allowed, but only when at least two positions are always opposite; in particular, a negatively weighted vertex cannot be a defect. The reason for this choice of weights is that a side attached by a bridge can contribute either a constant or a constant minus the colour at the attachment vertex. Section~\ref{sec:bridges} shows how to absorb the latter contribution into the weight at that vertex while preserving admissibility. This dependence between weights and boundary conditions is essential.

The \emph{spectrum} here means the finite set of attainable weighted sums $S(h)$ over admissible colourings, not a spectrum of eigenvalues. For an integer $L\geq0$, say that the model has \emph{spectral step at most $L$} if, for every pair of admissible colourings $h,g$, there is a finite sequence of admissible colourings
\begin{equation}\label{eq:walk}
 h=h_0,h_1,\ldots,h_k=g,\qquad
 \abs{S(h_{i+1})-S(h_i)}\leq L\quad(0\leq i<k).
\end{equation}
The colourings in a step may differ on arbitrarily many vertices. This is not a single-vertex reconfiguration condition.

\begin{lemma}\label{lem:potential}
Every boundary model admits an admissible colouring. For a nonempty finite spectrum of attainable sums, the spectral step bound $L$ is equivalent to the assertion that every gap between consecutive distinct attainable sums is at most $L$.
\end{lemma}
\begin{proof}
Let
\begin{equation}\label{eq:potential}
 M(h)=\sum_{uv\in E(H)}\ind{h(u)=h(v)}+\sum_v f_v(h(v)).
\end{equation}
If $h'$ is obtained by flipping $v$, precisely its $3-o_v$ non-opposite positions change their match status. Thus
\begin{equation}\label{eq:drop}
 M(h')-M(h)=3-o_v-2d_h(v).
\end{equation}
When $v$ is a defect this difference is strictly negative. A minimizer of the nonnegative integer $M$ is therefore admissible.

A sequence as in \eqref{eq:walk} connecting colourings of two consecutive distinct attainable sums must cross the gap between them in some step. Conversely, choose an admissible representative for every attainable sum between those of $h$ and $g$, order them by sum, and use $h,g$ as the endpoints. Equal sums can be joined in one step. This proves the equivalence, including singleton spectra.
\end{proof}

\begin{lemma}[One-defect repair]\label{lem:single}
Suppose that $h$ is admissible except possibly at a vertex $r$ with colour $s=h(r)$. There is an admissible colouring $g$ with
\begin{equation}\label{eq:one-delta}
 S(g)-S(h)\in\{0,-2s\}.
\end{equation}
If $r$ was a defect, $g(r)=-s$. If $h$ was already admissible, one may take $g=h$.
\end{lemma}
\begin{proof}
If $r$ is a defect, flip it. For any flipped defect $x$, its new monochromatic degree is $3-o_x-d_h(x)\leq1$. At a neighbouring vertex the degree decreases along an originally matching edge and increases along an originally opposite edge. There is at most one originally opposite position at a defect. Only if this position is an internal edge to a vertex $y$ can a new defect arise. If $y$ is a defect after the flip, flip $y$ next; otherwise stop. This procedure maintains at most one defect. Each flip strictly decreases \eqref{eq:potential}, so it terminates with an admissible colouring.

Let $x_1=r,x_2,\ldots,x_k$ be the flipped vertices in order and $s_i$ their colours just before the corresponding flips. The rule for the successor gives $s_i=(-1)^{i-1}s$. Whenever there is a successor, the two matching positions and the internal edge to that successor occupy all three positions. Consequently $o_{x_i}=m_{x_i}=0$ and $w_{x_i}=1$ for $i<k$. At the last vertex, at least two positions match, so $m_{x_k}\leq o_{x_k}\leq1$ and $w_{x_k}\in\{0,1\}$. The changes telescope to
\[
 S(g)-S(h)=-2s\left(\sum_{i=1}^{k-1}(-1)^{i-1}+w_{x_k}(-1)^{k-1}\right).
\]
The expression in parentheses equals $w_{x_k}$ when $k$ is odd and $1-w_{x_k}$ when $k$ is even. Both lie in $\{0,1\}$, proving \eqref{eq:one-delta}.

It remains to check the final colour of $r$. If there is no successor to $r$, it is flipped just once. Otherwise, immediately after the successor is flipped, $r$ has zero matching positions: its only match after its own flip was the edge to that successor. If it ever acquires a new match, consider the first flip causing one. Its degree becomes one. The flip causing that increase can propagate a defect only to $r$, so the chain stops at that flip. In particular $r$ cannot be flipped a second time. This proves the claimed final colour.
\end{proof}

\begin{corollary}[Changing one fixed input]\label{cor:input}
Replace one fixed boundary sign $t$ at a vertex by $-t$, starting from an admissible colouring. The new model has an admissible colouring with weighted sum changed by $0$ or $2t$.
\end{corollary}
\begin{proof}
The parameters $o_v,m_v,w_v$ do not change. The only possible defect is at the affected vertex and can arise only if its colour is the new sign $-t$. Apply Lemma~\ref{lem:single} with $s=-t$; if no defect arises, keep the original colouring.
\end{proof}

\section{Bridge contraction}\label{sec:bridges}
Consider one side $A$ of a cut bridge, with root $x$ at its endpoint. Replace the removed edge at $x$ by a single unmarked opposite position, so the budget \eqref{eq:budget} and weight at $x$ are unchanged. Reserve that position for an input $t\in\sgnset$; substituting $t$ makes the reserved position fixed with sign $t$. All other positions and weights remain unchanged. A colouring of the side is admissible at input $t$ if it is admissible in this substituted model.

\begin{lemma}[Profiles through a prescribed state]\label{lem:profiles}
For every admissible colouring $F_0$ of a side at an input $t_0$, there exist two admissible colourings $F(t)$, one for each input $t\in\sgnset$, with $F(t_0)=F_0$. They can be chosen to have one of the following forms, for a constant $C\in\Z$:
\begin{align}
 \text{fixed profile:}\quad &F(t)(x)=a, & S_A(F(t))&=C,\label{eq:fixed-profile}\\
 \text{opposite profile:}\quad &F(t)(x)=-t, & S_A(F(t))&=C-\delta t,\quad \delta\in\{0,1\}.\label{eq:opp-profile}
\end{align}
Here $a$ is a fixed sign, and $S_A$ uses the original weights of the side.
\end{lemma}
\begin{proof}
Keep $F(t_0)=F_0$. If the same colouring is admissible when the input is changed to $-t_0$, use it also for the second input. This is \eqref{eq:fixed-profile}. In this case the other two positions at the root have zero matches: for one of the two input signs, the input itself matches the root.

Otherwise, the only new defect is at the root, whose colour $s$ must be $-t_0$. Apply Lemma~\ref{lem:single} to the new boundary. The repaired root has colour $-s=t_0$, so the root is opposite to the input in both states. The two sums obey $S_A(F(s))-S_A(F(-s))\in\{0,-2s\}$. They therefore have the form \eqref{eq:opp-profile}: choose $\delta=0$ or $1$ according to the difference, and $C$ equal to the average of the two sums. This average is integral since the difference is even.
\end{proof}

A profile of $A$ can be substituted into its neighbouring side $B$ at the endpoint $y$ of the cut bridge. In \eqref{eq:fixed-profile}, replace the reserved position of $y$ by a fixed position of sign $a$. In \eqref{eq:opp-profile}, replace it by an opposite position with label $\delta$. In the latter case the weight at $y$ changes from $w_y$ to $w_y-\delta$; the position count remains three and the marked count remains at most the opposite count. In both cases the result is a boundary model of the same class.

For every admissible colouring $h_B$ of the contracted model, attach the chosen colouring $F(h_B(y))$ of $A$. At $x$ the actual bridge supplies precisely the input used to construct the profile. At $y$ the bridge has precisely the match status of the substituted boundary position. Thus admissibility is preserved at both ends and everywhere else. With $S_B^{\mathrm{new}}$ denoting the new weighted sum, the full colouring $h$ satisfies
\begin{equation}\label{eq:lift}
 S_{A\cup B}(h)=S_B^{\mathrm{new}}(h_B)+C.
\end{equation}
Indeed, in the opposite case the contribution from $A$ is $C-\delta h_B(y)$, exactly the term removed by the new weight; in the fixed case the contribution is constant. In particular, differences of sums are unchanged by this lifting operation.

\begin{lemma}[Bridge gluing]\label{lem:glue}
Suppose a connected boundary model is split by a bridge into sides $A$ and $B$. Assume that each side has spectral step at most $L$ for every replacement of its reserved bridge position by a fixed sign or an opposite position labelled $0$ or $1$. Then the original model has spectral step at most $L$.
\end{lemma}
\begin{proof}
Take arbitrary admissible colourings $f,g$ of the original model. Let $x\in A$ and $y\in B$ be the bridge endpoints. By Lemma~\ref{lem:profiles}, choose a profile of $A$ through $f|_A$ at input $f(y)$, and a profile of $B$ through $g|_B$ at input $g(x)$. Denote their root-output functions by $\alpha$ and $\beta$. Each function is constant or $t\mapsto-t$.

We claim there are signs $a,b$ with
\begin{equation}\label{eq:compatible}
 a=\alpha(b),\qquad b=\beta(a).
\end{equation}
If $\alpha$ is constant, take its value as $a$ and then let $b=\beta(a)$; if $\beta$ is constant, do the reverse. If both negate, take any opposite pair. Joining the selected rows with inputs $b$ on $A$ and $a$ on $B$ yields an admissible full colouring $k$.

Now contract $A$ using its profile through $f$. Both $f|_B$ and $k|_B$ are admissible in the contracted model: their lifted rows of $A$ are exactly $f|_A$ and $k|_A$. By hypothesis there is a sequence between those two restrictions with step at most $L$. Lift each member using the same profile. Equation~\eqref{eq:lift} shows that all sums change by exactly the differences of the contracted sums, giving the desired sequence from $f$ to $k$.

Similarly, contract $B$ using its profile through $g$. The restrictions $k|_A$ and $g|_A$ lift to $k$ and $g$, and the step bound is again preserved. Concatenation gives a sequence from $f$ to $g$ with step at most $L$. The two contractions need not use the same profile; the common state $k$ is what joins their lifted sequences.
\end{proof}

\section{A sharp spectral bound for cacti}\label{sec:cactus}
\begin{lemma}[Forests]\label{lem:forest}
Every forest boundary model has spectral step at most $4$. If, in addition,
\begin{equation}\label{eq:forest-extra}
 m_v+\deg_H(v)\leq2\qquad(v\in V(H)),
\end{equation}
then its spectral step is at most $2$.
\end{lemma}
\begin{proof}
For a one-vertex model, changing its colour changes the weighted sum by $2\abs{1-m_v}\leq4$. Under \eqref{eq:forest-extra} this is at most $2$. For a tree with an edge, split across any edge, which is a bridge. Every contraction on either side remains a smaller tree boundary model. Moreover, if \eqref{eq:forest-extra} holds, it is preserved: at the root the internal degree decreases by one and the marked count increases by at most one. Induction on the number of vertices, with Lemma~\ref{lem:glue}, proves both assertions for trees. For a forest, join two admissible colourings by changing the components one at a time, applying the respective tree bound. There are no internal edges between components and the weighted sum is additive. The empty forest is immediate.
\end{proof}

\begin{lemma}[Cycles]\label{lem:cycle}
A boundary model whose internal graph is a simple cycle has spectral step at most $4$.
\end{lemma}
\begin{proof}
Every cycle vertex has exactly one boundary position and weight either $0$ or $1$. Choose a vertex $v$ of weight zero if any exists; otherwise choose any vertex. Removing $v$ leaves a path $Q$. For a prescribed sign $s$ of $v$, replace the two edges to $v$ by two fixed positions of sign $s$ at their respective endpoints in $Q$. The endpoints are distinct because the cycle is simple and has length at least three. All remaining weights are unchanged.

An admissible colouring of $Q$ extends to a \emph{temporary} colouring of the cycle by assigning $s$ to $v$. Only $v$ can be a defect. Define a repair of each temporary colouring as follows. If it is admissible, do nothing. If $w_v=0$, the external position at $v$ is opposite, so a defect at $v$ has two matching cycle neighbours. Flipping $v$ removes both matches and cannot create a defect elsewhere; its weight is zero and the sum does not change. If $w_v=1$, use Lemma~\ref{lem:single}; its correction belongs to $\{0,-2s\}$. Denote the resulting admissible colouring by $R_s(T)$ for a temporary colouring $T$.

First connect two temporary colourings with the same prescribed $s$. Their path restrictions can be joined using Lemma~\ref{lem:forest}. If $w_v=0$, the path step is at most $4$, and the repair leaves each sum unchanged. If $w_v=1$, there is no zero weight anywhere on the original cycle by choice of $v$. Hence the path initially has no marked opposite positions, satisfies \eqref{eq:forest-extra}, and has step at most $2$. The corrections of two consecutive temporary states both belong to $\{0,-2s\}$ and differ by at most $2$. After repair, their sums differ by at most $2+2=4$. Thus the repaired images give a sequence with step at most $4$ in either case. In particular, two admissible full colourings agreeing at $v$ can be joined in this way.

Next let $f,g$ be admissible full colourings with $f(v)=s$ and $g(v)=-s$. Starting from $f|_Q$, change the two fixed positions representing $v$ from $s$ to $-s$, one at a time. After each change, use Corollary~\ref{cor:input} to obtain an admissible path colouring. Each change gives a correction in $\{0,2s\}$, even if the first repair changes the colour at the other endpoint. The intermediate path model may have different signs at its two new fixed positions; it need not correspond to a full cycle colouring.

After both changes, restore $v$ with sign $-s$. This gives a temporary state $T$ with
\[
 S(T)-S(f)=2js-2w_vs,\qquad j\in\{0,1,2\}.
\]
Repair $T$ using $R_{-s}$. If $w_v=0$, the correction is zero; if $w_v=1$, it is $0$ or $2s$. In either case
\[
 \frac{S(R_{-s}(T))-S(f)}{s}\in
 \begin{cases}
 \{0,2,4\}, & w_v=0,\\
 \{-2,0,2,4\}, & w_v=1.
 \end{cases}
\]
Hence $f$ and $R_{-s}(T)$ are joined by a step of size at most $4$. To continue, use the \emph{unrepaired} restriction $T|_Q$ as the starting path state and $g|_Q$ as the target, with both boundary inputs fixed at $-s$. Apply the same-sign argument to these temporary states, repairing each separately. The resulting full sequence starts at $R_{-s}(T)$ and ends at $g$. It is not necessary for the repaired colour at $v$ to remain $-s$. This distinction ensures the two parts of the construction have a common endpoint. The lemma follows.
\end{proof}

A cactus here means a graph in which any two simple cycles sharing an edge have the same edge set. For subcubic graphs this implies that distinct simple cycles are vertex-disjoint: two cycles through the same vertex use two incident edges each, and at degree at most three they must share an edge.

\begin{theorem}[Cactus spectrum]\label{thm:cactus}
Every connected subcubic cactus boundary model has spectral step at most $4$. The same bound holds for a disjoint union of such models. In particular, under arbitrary fixed boundary signs and unit weights the gaps between consecutive attainable ordinary imbalances are at most $4$. This constant is best possible.
\end{theorem}
\begin{proof}
Induct on the number of vertices in a connected cactus. If a bridge exists, each side of its cut is a smaller connected cactus, and its contracted models all satisfy \eqref{eq:budget}. The induction hypothesis and Lemma~\ref{lem:glue} give the bound $4$.

If there is no bridge and the graph has no edge, it is a single vertex. Otherwise every edge lies on a simple cycle: removing a non-bridge edge leaves a path between its endpoints, and a shortest such path together with the edge is a cycle. Choose one cycle $C$. Any edge at a vertex of $C$ lies on a cycle sharing that vertex with $C$. The degree bound forces the two cycles to share an edge, and the cactus property forces their edge sets to coincide. Hence every edge incident with a vertex of $C$ lies on $C$. Connectedness now implies that the graph is exactly $C$. The one-vertex case follows from Lemma~\ref{lem:forest}, and the cycle case from Lemma~\ref{lem:cycle}. This proves the induction without requiring any restriction on the parity of the cycles. For disconnected graphs, concatenate the sequences on the connected components, using additivity as in Lemma~\ref{lem:forest}.

For sharpness, take the star $K_{1,3}$ and give each leaf one fixed positive and one fixed negative position. All weights are one. Each leaf already has one matching fixed position, so its colour must be opposite to that of the centre. Both choices of centre colour are admissible and the only attainable sums are $-2$ and $+2$. Their gap is $4$.
\end{proof}

\section{Reduction by even-cycle packing}\label{sec:packing}
\begin{lemma}\label{lem:packing}
Every finite simple graph admits a family of pairwise vertex-disjoint even cycles such that the induced graph on the uncovered vertices has no even cycle and is a cactus (possibly disconnected).
\end{lemma}
\begin{proof}
Choose a family of pairwise vertex-disjoint simple even cycles that is maximal under inclusion. Such a family exists by finiteness, and the empty family is permitted. Let $U$ be the union of its vertex sets and $W=V(G)\setminus U$. An even cycle in $G[W]$ could be added to the family, a contradiction.

To see the cactus property, suppose two simple cycles $C,C'$ in $G[W]$ share an edge but have different edge sets. Both cycles are odd. Since their edge sets differ, $C'$ contains an edge not in $C$: a simple cycle cannot be a proper subgraph of another simple cycle. Traversing $C'$, there is a segment outside the edges of $C$ from a vertex $a$ of $C$ to another vertex $b$ of $C$, with no internal vertex on $C$. The endpoints are distinct because $C'$ is simple and contains an edge of $C$. This segment may be a single chord. The segment together with the two $a$--$b$ arcs of $C$ gives three internally vertex-disjoint paths between distinct endpoints. Two of their lengths have the same parity. Their union is a simple even cycle (the paths are distinct and cannot both be the same single edge in a simple graph), a contradiction. Thus the two cycles have identical edge sets, as required.
\end{proof}

For the rest of the proof let $G$ be a subcubic graph of even order. At each vertex with degree less than three, introduce $3-\deg_G(v)$ unmarked opposite positions. All weights are one. This padding changes neither the monochromatic degree in \eqref{eq:ordinary} nor the ordinary sum. When $G$ is cubic, no padding is necessary.

Choose a packing as in Lemma~\ref{lem:packing}. For each fixed colouring of $U$, the residue $G[W]$ becomes a boundary model: replace every edge from a vertex of $W$ to a vertex of $U$ by one fixed position bearing the colour of its endpoint in $U$, and retain all the original opposite positions. Distinct edges give distinct positions even when several such edges meet one vertex. All weights in the residue are still one. For every $v\in W$, its monochromatic degree in the boundary model equals its monochromatic degree in the full graph. Write
\begin{equation}\label{eq:residue}
 P(h)=\sum_{v\in W}h(v).
\end{equation}
By Lemma~\ref{lem:potential} the residue has an admissible colouring for every colouring of $U$, and by Theorem~\ref{thm:cactus} it has spectral step at most $4$ for each fixed boundary. Here $P$ is always even: both $\abs{V(G)}$ and $\abs U$ are even, so $\abs W$ is even, and each summand is congruent to one modulo two.

\section{Repair at two marked edges}\label{sec:marked}
We first prove a repair lemma allowing several defects. This will be applied to the full graph with unit weights.

\begin{lemma}[A repair chain among defects]\label{lem:multiple}
Let $H$ be a boundary model all of whose weights are one, and let $h$ be any colouring. If $r$ is a defect of colour $s$, there is a colouring $g$ such that
\begin{equation}\label{eq:multi-delta}
 S(g)-S(h)\in\{0,-2s\},
\end{equation}
and every defect $v$ of $g$ satisfies
\begin{equation}\label{eq:survivor}
 v\ne r,\qquad g(v)=h(v),\qquad
 d_g(v)+\ind{rv\in E(H),\ h(r)=h(v)}\leq d_h(v).
\end{equation}
In particular, no new defect survives, and the degree of a surviving defect originally joined to $r$ by a monochromatic edge decreases by at least one.
\end{lemma}
\begin{proof}
Start by flipping $r$. At each step follow the unique formerly opposite internal edge if its other endpoint is a defect after the flip; otherwise stop. The successor must be processed immediately even if it was already a defect before the current flip. This is the same rule as in Lemma~\ref{lem:single}, now allowing other defects to be present. By \eqref{eq:drop} the potential strictly decreases at every flip, so termination is guaranteed. The colours at the flipped vertices, just before their respective flips, alternate beginning with $s$. Since all weights are one, the sum of all changes is $-2s$ for an odd number of flips and zero for an even number. This gives \eqref{eq:multi-delta}.

For clarity, track a set of inactive defects. Initially it consists of all original defects except $r$, which is the active defect. Immediately after flipping an active vertex, its own degree is at most one. All neighbours on previously matching edges lose one match. There is at most one neighbouring vertex gaining a match; if that vertex has degree at least two, it becomes the next active defect. Delete from the inactive set every vertex which is now admissible, and delete the next active vertex if there is one. Inductively, the inactive set at the start of each step is exactly the set of current defects other than the active vertex. Each member retains its original colour, has never been flipped, and its monochromatic degree has not increased since the beginning of the chain. No vertex is added to this set. Indeed, the only possible newly created defect is immediately assigned to be active. The final set therefore consists of a subset of the original inactive defects. If one of these was originally joined to $r$ by a monochromatic edge, its degree decreased at the first flip and never increased while it remained inactive. This proves \eqref{eq:survivor}. No assertion about the colours of vertices which are already admissible at termination is needed.
\end{proof}

For an edge $e=ab$, let $I_e(v)=\ind{v=a\text{ or }v=b}$. Two different edges may have a common endpoint.

\begin{lemma}[Clearing one marked edge]\label{lem:clear}
In a unit-weight boundary model, let $e=ab$ be an internal edge and let $B:V(H)\to\N$ have $B(v)\geq1$ everywhere. Suppose that a colouring $h$ satisfies
\begin{equation}\label{eq:clear-cap}
 d_h(v)\leq B(v)+I_e(v)\qquad(v\in V(H)),
\end{equation}
and all defective endpoints of $e$ have the same prescribed colour $s$. Then there is a colouring $g$ with $d_g(v)\leq B(v)$ everywhere and $S(g)-S(h)\in\{0,-2s\}$. Every defect of $g$ is an original defect of $h$ with its original colour.
\end{lemma}
\begin{proof}
If $d_h(v)\leq B(v)$ for all $v$, take $g=h$. Otherwise, by \eqref{eq:clear-cap}, some endpoint $r$ of $e$ satisfies $d_h(r)>B(r)\geq1$. It is a defect of colour $s$. Apply Lemma~\ref{lem:multiple} starting at $r$. If some $v$ still satisfies $d_g(v)>B(v)$, then $v$ is a surviving defect. By \eqref{eq:survivor}, $v\ne r$ and $d_h(v)\geq d_g(v)>B(v)$. Equation~\eqref{eq:clear-cap} forces $v$ to be the other endpoint of $e$. Both endpoints were defective and of colour $s$ in $h$, so the indicator in \eqref{eq:survivor} equals one. It follows that $d_g(v)+1\leq d_h(v)\leq B(v)+1$, a contradiction. The sum and the assertion about surviving defects follow from Lemma~\ref{lem:multiple}.
\end{proof}

\begin{lemma}[Two marked edges]\label{lem:two}
Let $e_1,e_2$ be internal edges of a unit-weight boundary model. Suppose a colouring $h$ has both endpoints of $e_1$ of colour $s$ and both endpoints of $e_2$ of colour $t$, and that
\begin{equation}\label{eq:two-hyp}
 d_h(v)\leq1+I_{e_1}(v)+I_{e_2}(v)\quad(v\in V(H)),
 \qquad S(h)=s+t.
\end{equation}
Then there is an admissible colouring $g$ with $\abs{S(g)}\leq2$. The edges are allowed to share an endpoint.
\end{lemma}
\begin{proof}
Apply Lemma~\ref{lem:clear} to $e_1$ with $B(v)=1+I_{e_2}(v)$. Obtain a colouring $k$ with degrees at most this budget and sum $S(k)=S(h)-2\varepsilon_1 s$ for some $\varepsilon_1\in\{0,1\}$. Every surviving defective endpoint of $e_2$ has its original colour $t$. Thus Lemma~\ref{lem:clear} applies again to $e_2$, now with the budget identically one. We obtain an admissible colouring $g$ and $\varepsilon_2\in\{0,1\}$ with
\[
 S(g)=s+t-2\varepsilon_1s-2\varepsilon_2t
     =(1-2\varepsilon_1)s+(1-2\varepsilon_2)t.
\]
Each summand is a sign, giving the desired bound. The endpoint indicators and the two successive budgets account for shared endpoints, including an initial defect of degree three. No disjointness of the repair chains or of the two marked edges is used.
\end{proof}

\section{Proof of the main theorem}\label{sec:mainproof}
Use the packing and residue from Section~\ref{sec:packing}. On each packed cycle choose a cyclic ordering $v_0,\ldots,v_{n-1}$ with $n$ even and assign the alternating signs $(-1)^i$ at $v_i$. Every vertex of $U$ has two cycle neighbours of the opposite colour and is therefore admissible, regardless of its remaining incidence. Extend this colouring to $W$ using Lemma~\ref{lem:potential}. Call the resulting full admissible colouring $f$ and put $P_0=P(f)$.

\subsection{A partially inverted cycle}
The following observation identifies the precise hypotheses of Lemma~\ref{lem:two}. Consider a state in which $W$ is admissible at the current boundary, every packed cycle except possibly one is alternating, and the remaining cycle is obtained from its original alternating colouring by inverting the first $k$ vertices in cyclic order. We allow $k=0$ and $k=n$; in those cases all cycles are alternating.

Suppose first that $P=0$ and all packed cycles are alternating. The full sum is zero and the colouring is admissible. Otherwise, if $P=0$ and $0<k<n$, the two marked edges are
\[
 e_1=v_{n-1}v_0,\qquad e_2=v_{k-1}v_k.
\]
Both are monochromatic: exactly one of their endpoints has been inverted. All other edges belonging to the chosen cycles have endpoints of opposite colours. Let $s$ and $t$ be the colours of $e_1$ and $e_2$, respectively. Then $t=(-1)^{k-1}s$. The original first sign at $v_0$ was $-s$. Consequently, the sum on the affected cycle after the inversion equals
\begin{equation}\label{eq:prefix}
 -2\sum_{i=0}^{k-1}(-s)(-1)^i
 =\begin{cases}0,&k\text{ even},\\2s,&k\text{ odd},\end{cases}
 =s+t.
\end{equation}
The other packed cycles contribute zero and $P=0$, so the full sum is $s+t$.

At each vertex of $U$, there is at most one incidence outside its two chosen cycle edges. The two exceptional edges are the only possible monochromatic chosen cycle edges. At vertices of $W$, admissibility holds with respect to the actual colours of $U$. Hence at every vertex
\[
 d_h(v)\leq1+I_{e_1}(v)+I_{e_2}(v).
\]
This remains true when there are chords, edges between distinct packed cycles, or shared endpoints of $e_1$ and $e_2$. Thus \eqref{eq:two-hyp} holds for the full unit-weight boundary model. Lemma~\ref{lem:two} gives a 2-quasi-bisection. In this last repair, the previously maintained value of $P$ and the alternating pattern on other cycles need not be preserved.

\subsection{The passage through the residue spectrum}
Starting at $f$, invert the packed cycles one at a time, flipping their vertices consecutively in cyclic order. After each vertex flip, restore admissibility on $W$ before proceeding. At a flipped vertex of $U$, two incidences belong to its chosen cycle; hence at most one edge from it can go to $W$. Accordingly, the boundary model of $W$ changes in at most one fixed position. Corollary~\ref{cor:input} guarantees a restored admissible residue colouring with change of $P$ at most two in absolute value. If no position changes, keep the residue colouring. During this stage defects on $U$ are permitted; the only invariant concerning admissibility is on $W$.

If $P=0$ at any state of this stage (including its initial or final state), the preceding argument gives the required full colouring. Otherwise, because $P$ is even and changes by at most two at a step, its sign cannot change during this stage. In particular $P_0\ne0$ and the final residue sum has the same sign as $P_0$.

After all packed cycles have been inverted, fix their new signs. They again alternate and contribute zero to the full sum. The colouring $-f|_W$ is admissible at this final boundary: relative to $f$, both endpoints of every internal edge of $W$ and both signs compared at every fixed position have been inverted. The padded opposite positions remain opposite. Its sum is $-P_0$.

Apply Theorem~\ref{thm:cactus} to the residue at this fixed final boundary. It joins the attained residue colouring to $-f|_W$ through admissible colourings in which successive values of $P$ differ by at most four. Every full colouring in this sequence is admissible because the packed cycles now alternate. The sequence begins and ends with sums of opposite signs. If every sum had absolute value greater than two, then all sums, being even, would be at least $4$ or at most $-4$. A step changing sign would have size at least eight, contradicting the bound four. Thus some state has $\abs P\leq2$. Its full sum is precisely $P$, proving Theorem~\ref{thm:subcubic}.

The argument also covers an empty packing (no inversion steps) and an empty residue (so $P=0$), including the empty graph. No step requires $G$ to be connected. For a cubic graph the order is even by the handshaking identity $3\abs{V(G)}=2\abs{E(G)}$. Theorem~\ref{thm:main} follows.

\begin{corollary}\label{cor:multi}
Every finite cubic multigraph, with loops counted twice in the degree if permitted, admits a partition of its vertices into two classes of sizes differing by at most two, with each monochromatic component having at most two vertices.
\end{corollary}
\begin{proof}
The order is even by the handshaking identity. Delete loops and replace each nonempty class of parallel edges by one edge to obtain a simple subcubic graph on the same vertices. Apply Theorem~\ref{thm:subcubic}. For each colour class, deleting loops and identifying parallel edges does not change which vertices lie in the same connected component. Thus the same partition satisfies the stated component-order bound in the multigraph.
\end{proof}
\begin{remark}
Corollary~\ref{cor:multi} concerns component order, not monochromatic incidence degree or acyclicity. A monochromatic loop or parallel edges may remain in a component of order at most two. The stronger incidence-degree formulation is used only for simple graphs.
\end{remark}

\section{Correspondence with the Lean formalization}\label{sec:lean}
The accompanying project uses Lean~4.36.0 and Mathlib~v4.36.0. In the namespace \lean{QuasiBisection}, the declarations \lean{subcubic_quasi_bisection} and \lean{subcubic_quasi_bisection_classes} establish Theorem~\ref{thm:subcubic}, respectively in terms of Boolean colours and their imbalance, and in terms of colour-class cardinalities and monochromatic component orders. The cubic specializations \lean{cubic_quasi_bisection} and \lean{cubic_quasi_bisection_classes} give Theorem~\ref{thm:main} by the handshaking identity. Table~\ref{tab:lean} identifies the principal corresponding declarations under \lean{formal/QuasiBisection/}.

In \lean{BoundaryModel}, the fields \lean{fixed}, \lean{opposite} and \lean{marked} implement the parameters of \eqref{eq:budget}, and \lean{Bridge.Walk} expresses the bounded-step chains of \eqref{eq:walk}. The Lean proof of the main theorem uses a contradiction equivalent to the passage through zero or a small imbalance in Section~\ref{sec:mainproof}. Lean also verifies the upper bound on the cactus step. The multigraph corollary and the sharpness example are proved in the text.

\begin{table}[htbp]
\centering\small
\begin{tabularx}{\textwidth}{@{}>{\raggedright\arraybackslash}p{.23\textwidth}>{\raggedright\arraybackslash}p{.26\textwidth}>{\raggedright\arraybackslash}X@{}}
\toprule
Mathematical result & Lean module & Principal declaration(s) \\
\midrule
Lemma~\ref{lem:potential} (existence) and Lemma~\ref{lem:single} & \lean{DefectRepair}, \lean{SingleDefect} & \lean{repair_exists}, \lean{single_repair_trace}, \lean{single_repair} \\
Lemma~\ref{lem:profiles} & \lean{BoundaryProfiles} & \lean{boundaryProfiles} \\
Lemma~\ref{lem:glue} & \lean{BridgeGluing}, \lean{BridgeContraction}, \lean{BridgeCut} & \lean{bridgeGluing}, \lean{boundaryBridge}, \lean{BridgeCut.spectralStep} \\
Lemma~\ref{lem:forest} & \lean{ForestSpectrum} & \lean{forest_step_four}, \lean{forest_step_two} \\
Lemma~\ref{lem:cycle} & \lean{CycleSpectrum}, \lean{CycleGraph} & \lean{open_cycle_step}, \lean{cycle_step_four} \\
Theorem~\ref{thm:cactus} (step bound) & \lean{CactusSpectrum}, \lean{ComponentSpectrum} & \lean{cactus_step_four}, \lean{cactus_forest_step_four} \\
Lemma~\ref{lem:packing} & \lean{EvenCyclePacking}, \lean{OddCactus} & \lean{exists_even_packing_cactus} \\
Residue and boundary change & \lean{ResidualBoundary}, \lean{PackingBoundary} & \lean{residue_step_four}, \lean{packing_change_exterior} \\
Lemmas~\ref{lem:multiple}--\ref{lem:two} & \lean{SingleDefect}, \lean{MarkedEdges} & \lean{repair_among_defects}, \lean{repair_two_edges} \\
Partial inversion and completion & \lean{CycleOrder}, \lean{InversionRepair}, \lean{MainTheorem} & \lean{prefix_imbalance}, \lean{finish_partial}, \lean{finish_settled}, \lean{subcubic_quasi_bisection} \\
\bottomrule
\end{tabularx}
\caption{Correspondence between the mathematical proof and the Lean development.}
\label{tab:lean}
\end{table}

\paragraph{Reproducibility.} The manifest pins Mathlib to revision \lean{fc4f59e577feb1a4956671bc04ac18aaaf8cd842}. The ancillary files contain the Lean sources and pinned manifest. From the \lean{formal/} directory in the ancillary material, run
\begin{verbatim}
lake exe cache get
lake build
lake env lean -DwarningAsError=true Smoke.lean
export LEAN_NUM_THREADS=1
lake env leanchecker --fresh -v QuasiBisection.MainTheorem
\end{verbatim}
The recorded \texttt{\#print axioms} output for both formulations of the subcubic theorem and their cubic specializations lists only \lean{propext}, \lean{Classical.choice}, and \lean{Quot.sound}. The analytic proof is independent of finite computational experiments.

\paragraph{Use of AI tools.}
We used multiple AI agents, including Kimi and Codex, to explore proof
strategies, develop and review mathematical arguments, assist with the
Lean formalization, and prepare the manuscript. After each AI-confirmed
step, we manually verified the corresponding lemma, proof, or theorem
before proceeding. We directed the research by updating the agents'
instructions as the work progressed. Claude Opus additionally reviewed
the formalization.
We verified that the formal statement matches the theorem in the paper,
and Lean's kernel checked the formal proof. The accompanying AI workflow
note describes the agents' roles, discarded approaches, and verification
steps.

\end{document}